\documentclass{article}
\usepackage{amsfonts}
\usepackage{amsmath}

\input{tcilatex}
\begin{document}

\begin{center}
HARNACK, H\"{O}LDER, GAUSS AND WIDDER: SERRIN'S PARABOLIC LEGACY
\end{center}

\bigskip

D. G. Aronson

School of Mathematics

University of Minnesota

Minneapolis MN 55455

e-mail: arons001@umn.edu

\bigskip

ABSTRACT: James Serrin's fundamental contributions to the theory of
quasilinear elliptic equations are well-known and widely appreciated. He
also made less well-known contributions to the theory of quasilinear
parabolic equations which we discuss in this note. J\"{u}gen Moser gave
greatly simplified proofs of the De Giorgi-Nash regularity results for
linear divergence structure elliptic and parabolic differential equations
using an original iterative technique. Serrin extended Moser's techniques
and applied them to the study of divergence structure quasilinear elliptic
and, in collaboration with Aronson, to divergence structure quasilinear
parabolic equations. Specifically, among other results, they proved a
maximum principle, H\"{o}lder continuity of generalized solutions and
derived a Harnack principle for a very broad class of quasilinear parabolic
equations. In subsequent work, Aronson applied these results to study
non-negative solutions to divergence structure linear equations without
regularity assumptions on the coefficients. The results include a two-sided
Gaussian estimate for the fundamental solution and a generalization of the
Widder Representation Theorem.

\bigskip

James Serrin's fundamental contributions to the theory of quasilinear
elliptic equations [14] are well-known and widely appreciated. He also made
less well-known but important contributions to the theory of quasilinear
parabolic equations which we will describe here.

\bigskip

In two remarkable, essentially simultaneous works, Ennio De Giorgi [6]
proved in 1956-57 the H\"{o}lder continuity of weak solutions of the
divergence structure elliptic equation%
\begin{equation*}
\left\{ A_{ij}(x)u_{x_{i}}\right\} _{x_{j}}=0
\end{equation*}%
and John Nash [11] proved in 1957-58 the same result for weak solutions of
the parabolic equation%
\begin{equation*}
u_{t}-\left\{ A_{ij}(x,t)u_{x_{i}}\right\} _{x_{j}}=0.
\end{equation*}%
Here and throughout this note we employ the convention of summation over
repeated indicies. In case everything is independent of $t$, Nash also
derives the De Giorgi result for elliptic equations. It should be noted that
aside from boundedness, measurability and uniform ellipticity or
parabolicity there are no further assumptions on the coefficients. Their
work was totally independent and their methods completely different.
Subsequently J\"{u}rgen Moser [9], [10] introduced an iterative techniques
to prove Harnack inequalities for both the elliptic (1961) and parabolic
equations (1964) which enabled him to give greatly simplified proofs of the
De Giorgi-Nash continuity results. In 1964 James Serrin [14] extended
Moser's techniques to obtain a Harnack inequality and H\"{o}lder continuity
for a broad class of divergence structure quasilinear elliptic equations of
the form%
\begin{equation*}
\func{div}\mathcal{A(}x,u,u_{x})+\mathcal{B}(x,u,u_{x})=0.
\end{equation*}%
We will be concerned with the extensions of Serrin's work to parabolic
equations and its ramifications.

Let $\Omega $ be a bounded domain in $\mathbf{R}^{n}$ and consider the
space-time cylinder $Q=\Omega \times (0,T)$ for some fixed $T>0$. We treat
the second order quasilinear equation%
\begin{equation}
u_{t}=\func{div}\mathcal{A(}x,t,u,u_{x})+\mathcal{B}(x,t,u,u_{x})\text{,} 
\tag{1}
\end{equation}%
where $\mathcal{A=(A}_{1},...,A_{n})$ is a given vector function of $%
(x,t,u,u_{x})$, $\mathcal{B}$ is a given scalar function of the same
variables, and $u_{x}=(\partial u/\partial u_{x_{1}},...\partial u/\partial
u_{x_{n}})$ denotes the spatial gradient of the dependent variable $u=u(x,t)$%
. Also here $\func{div}\mathcal{A}$ refers to the divergence of the vector $%
\mathcal{A}(x,t,u(x,t),u_{x}(x,t))$ with respect to the variables $%
(x_{1},...,x_{n})$. The structure of (1) is determined by the functions $%
\mathcal{A}(x,t,u,p)$ and $\mathcal{B}(x,t,u,p)$. We assume that they are
defined and measurable for all $(x,t)\in Q$ and for all values of $u$ and $p$%
. In addition, $\mathcal{A}$ and $\mathcal{B}$ will be required to satisfy
certain inequalities whose description follows.

A function $w=w(x,t)$ which is defined and measurable on $Q$ belongs to the
Bochner space $L^{p,q}(Q)$ if the iterated integral%
\begin{equation*}
\left\Vert w\right\Vert _{p,q}=\left\{ \int_{0}^{T}(\tint_{\Omega
}\left\vert w\right\vert ^{p}dx)^{q/p}dt\right\} ^{1/q}
\end{equation*}%
is finite. Here $p,q$ may be any real numbers $\geq 1$; and with the obvious
use of $L^{\infty }$ norms rather than integrals we can allow $p$ and $q$ to
have the value $\infty $.

We shall always assume that the functions $\mathcal{A}$ and $\mathcal{B}$
satisfy inequalities of the form%
\begin{eqnarray}
p\cdot \mathcal{A}(x,t,u,p) &\geq &a\left\vert p\right\vert
^{2}-b^{2}u^{2}-f^{2}  \notag \\
\left\vert \mathcal{B(}x,t,u,p)\right\vert &\leq &c\left\vert p\right\vert
+d\left\vert u\right\vert +g  \TCItag{2} \\
\left\vert \mathcal{A}(x,t,u,p)\right\vert &\leq &\overline{a}\left\vert
p\right\vert +e\left\vert u\right\vert +h.  \notag
\end{eqnarray}%
Here $a$ and $\overline{a}$ are positive constants, while the coefficients $%
b,c,...,h$ are non-negative functions of $(x,t)$ each contained in some
space $L^{p,q}(Q)$, where $p$ and $q$ are non-negative real numbers
(possibly\ different from coefficient to coefficient) such that 
\begin{equation}
p>2\text{ and }\frac{n}{2p}+\frac{1}{q}<\frac{1}{2}\text{ for }b,c,e,f,h 
\tag{3}
\end{equation}%
and%
\begin{equation}
p>1\text{ and }\frac{n}{2p}+\frac{1}{q}<1\text{ for }d,g\text{.}  \tag{4}
\end{equation}%
The norms of the functions $b,c,...,h$ in their respective spaces will be
denoted simply by $\left\Vert b\right\Vert ,\left\Vert c\right\Vert
,...,\left\Vert h\right\Vert $. It is clear that there exists a positive
number $\theta $ such that%
\begin{equation*}
p\geq \frac{2}{1-\theta }\text{ and }\frac{n}{2p}+\frac{1}{q}\leq \frac{%
1-\theta }{2}\text{ for }b,c,e,f,h
\end{equation*}%
and%
\begin{equation*}
p\geq \frac{1}{1-\theta }\text{ and }\frac{n}{2p}+\frac{1}{q}\leq 1-\theta 
\text{ for }d,g.
\end{equation*}%
The qualitative structure of equation (1) is entirely determined by the
coefficients $a$ and $\overline{a}$, the value of the norms $\left\Vert
b\right\Vert $ through $\left\Vert h\right\Vert $ in their respective spaces
and the numbers $\theta $ and $n$. A constant will be said to depend on the
structure of equation (1) if it is determined by these quantities (and is
uniformly bounded whenever these quantities are.)

The linear equation%
\begin{equation}
u_{t}=\left\{ A_{ij}(x,t)u_{x_{i}}+A_{j}(x,t)u+F_{j}(x,t)\right\}
_{x_{j}}+B_{j}(x,t)u_{x_{j}}+C(x,t)u+G(x,t)  \tag{5}
\end{equation}%
satisfies the preceding hypotheses if (i) the $A_{ij}$ are bounded and
measurable on $Q$ and there exists a constant $\nu >0$ such that $%
A_{ij}(x,t)\xi _{i}\xi _{j}\geq \nu \left\vert \xi \right\vert ^{2}$ almost
everywhere in $Q$ and for all $\xi \in \mathbf{R}^{n}$; (ii) the
coefficients $A_{j},B_{j}$ and $F_{j}$ each belong to some space $L^{p,q}(Q)$
with $p,q$ satisfying (3); and (iii) the coefficients $C$ and $G$ each
belong to some space $L^{p,q}(Q)$ with $p,q$ satisfying (4).

Without further hypotheses on $\mathcal{A}$ and $\mathcal{B}$ it is not
possible, in general, to speak of a classical solution of equation (1), and
it is correspondingly necessary to introduce the notion of generalized
solution. Let $u=u(x,t)$ be a function which is locally of class $%
L^{2,\infty }$ in $Q$ and possesses a strong derivative $u_{x}$ which is
locally of class $L^{2,2}$ in $Q$. Then $u$ will be called a \textit{weak
solution} of equation (1) in $Q$ if%
\begin{equation}
\tiint_{Q}\left\{ -u\varphi _{t}+\varphi _{x}\cdot \mathcal{A(}%
x,t,u,u_{x})\right\} dxdt=\tiint_{Q}\varphi \mathcal{B(}x,t,u,u_{x})dxdt 
\tag{6}
\end{equation}%
for any continuously differentiable function $\varphi =\varphi (x,t)$ having
compact support in $Q$. It would be simpler to consider a less general class
of weak solutions, namely those for which both $u_{x}$ and $u_{t}$ are
locally of class $L^{2,2}$ in $Q.$ However, from the point of view of
existence theory the latter class is not as natural as the former. Thus even
in the linear case some smoothness of the coefficients must be assumed in
order to prove the existence of solutions having strong time derivatives.

In [5] we prove a global maximum principle, a local boundedness theorem and
a Harnack inequality for weak solutions of equation (1) under the hypotheses
listed above. All of these results depend on certain integral inequalities
which are derived from the weak form of the differential equation. Moreover,
in each case the proof depends on the iterative techniques introduced by
Moser [9], [10] and further developed by Serrin [14] and Aronson \& Serrin
[4], [5]. In [4] iterative arguments are employed to prove a maximum
principle for equations whose structure inequalities are somewhat more
general than those specified in (2). In [5] we also prove the interior H\"{o}%
lder continuity of weak solutions and to study the growth properties of
non-negative solutions near certain parts of the boundary of their domain of
definition.

As in [14], the fundamental inequalities are derived by using powers of $u$
as test functions in the weak form (6) of the differential equation (1). For
parabolic equations this presents certain technical difficulties which stem
from the lack of an \textit{a priori} bound on $u$ and the lack of
assumptions on the time derivative of $u$. These difficulties are resolved
in great detail in [4] and [5]. Let $\eta =\eta (x,t)$ be a piecewise smooth
non-negative function which vanishes in a neighborhood of the parabolic
boundary $\Gamma $ of $Q$. Then for almost all values of $\tau \in (0,T)$%
\begin{eqnarray}
&&\frac{1}{\beta +1}\tint_{\Omega }\eta ^{2}\left\{ \overline{u}^{\beta
-1}-(b+1)\kappa ^{\beta }\overline{u}+\beta \kappa ^{\beta +1}\right\}
_{t=\tau }dx+\frac{a\beta }{2}\tiint \eta ^{2}\overline{u}^{\beta -1}|%
\overline{u}_{x}|^{2}dxdt\leq  \notag \\
&&\text{ \ \ \ \ \ \ \ \ \ \ \ \ \ \ \ \ \ \ \ \ \ \ \ \ \ \ }\tiint 
\mathfrak{F}\overline{u}^{\beta +1}dxdt+\frac{2}{\beta +1}\tiint \eta
\left\vert \eta _{t}\right\vert \overline{u}^{\beta +1}dxdt  \TCItag{7}
\end{eqnarray}%
where $\beta \geq 1$ and $\kappa >0$ are fixed constants and $\overline{u}%
=\max (0,u)+\kappa $. The integrations extend over $\Omega \times (0,\tau ).$
Note that $\kappa =0$ if $f,g,h$ are zero. Moreover%
\begin{equation*}
\mathfrak{F}=F\eta ^{2}+2G\eta \left\vert \eta _{x}\right\vert +H\left\vert
\eta _{x}\right\vert ^{2}
\end{equation*}%
with%
\begin{equation*}
F=\beta \left( b^{2}+\frac{f^{2}}{\kappa ^{2}}\right) +\left( d+\frac{g}{%
\kappa }\right) +\frac{c^{2}}{a},G=e\frac{h}{\kappa }\text{ and }H=\frac{%
4a^{\prime 2}}{a}.
\end{equation*}%
This inequality and its variants form the basis of the iterative arguments
which are employed in [4] and [5]. The set $\Gamma =\{\partial \Omega \times
\lbrack 0,T)\}\cup \{\Omega $ $\times (t=0)\}$ is called the parabolic
boundary of $Q$ and we say that $u\leq M$ on $\Gamma $ if for every $%
\varepsilon >0$ there is a neighborhood of $\Gamma $ in which $u\leq
M+\varepsilon $.

\textit{MAXIMUM PRINCIPLE:\ Let }$u$\textit{\ be a weak solution of (1) in }$%
Q$\textit{\ such that }$u\leq M$\textit{\ on }$\Gamma $\textit{. Then almost
everywhere in }$Q$%
\begin{equation*}
u(x,t)\leq M+Ck,
\end{equation*}%
\textit{where }%
\begin{equation*}
k=(\left\Vert b\right\Vert +\left\Vert d\right\Vert )\left\vert M\right\vert
+(\left\Vert f\right\Vert +\left\Vert g\right\Vert )
\end{equation*}%
\textit{and }$C$\textit{\ depends only on }$T,\left\vert \Omega \right\vert $%
\textit{\ and the structure of (1). }

Actually\textit{\ \ \ \ \ \ \ \ \ \ \ \ \ \ \ \ \ \ \ \ \ \ \ \ \ \ \ \ \ \
\ \ \ \ \ \ \ \ \ \ \ \ \ \ \ \ \ \ \ \ \ \ \ \ \ \ \ \ \ \ \ \ } 
\begin{equation*}
C=C(a,\left\Vert b\right\Vert ,\left\Vert c\right\Vert ,\left\Vert
d\right\Vert ,\theta ,n,T,\left\vert \Omega \right\vert )
\end{equation*}%
and so, in particular, is independent of $a^{\prime },e,f,g$ \ and $h.$ Also
the structural inequality $\left\vert \mathcal{A}(x,t,u,p)\right\vert \leq
a^{\prime }\left\vert p\right\vert +e\left\vert u\right\vert +h$ does not
enter into the proof. All that is really needed is that the term $\varphi
_{x}\cdot \mathcal{A}$ be integrable for $\varphi _{x}$ belonging to $%
L^{2,2} $ locally in $Q$. Finally note that the maximum principle remains
true if the differential equation (1) is replaced by the differential
inequality%
\begin{equation*}
u_{t}\leq \func{div}\mathcal{A(}x,t,u,u_{x})+\mathcal{B}(x,t,u,u_{x}).
\end{equation*}%
Under the same hypothesis as in the maximum principle there is also a
minimum principle. Specifically, if $u\geq M$ on $\Gamma $ then%
\begin{equation*}
u(x,t)\geq M-Ck
\end{equation*}%
in $Q$.

The proof of the maximum principle involves iteration procedures based on
the fundamental inequality (7) and one of its variants. Let $r=(\beta +1)/2$
and $v=\overset{-}{u}^{r}$. Introduce the norm%
\begin{equation*}
|||v|||=\sup \left\Vert v\right\Vert _{p^{\prime },q^{\prime }}
\end{equation*}%
where the supremum is taken over all exponent pairs $(p^{\prime },q^{\prime
})$ whose H\"{o}lder conjugates $(p,q)$ satisfy%
\begin{equation*}
\frac{n}{2p}+\frac{1}{q}\leq 1-\theta \text{ and }p\geq \frac{1}{1-\theta }.
\end{equation*}%
Using (7), it is shown that%
\begin{equation*}
|||v^{\sigma }|||^{2/\sigma }\leq Cr^{2}|||v|||^{2}
\end{equation*}%
where $\sigma =1+2\theta /n$ and $C$ depends only on the structure of (1).
For $m=0,1,2,...,$ set 
\begin{equation*}
\varphi _{m}=|||v|||^{2/r}
\end{equation*}%
and rewrite the preceding inequality as%
\begin{equation*}
\varphi _{m+1}\leq C^{1/\sigma ^{m}}\sigma ^{2m/\sigma ^{m}}\varphi _{m}.
\end{equation*}%
Iteration yields%
\begin{equation*}
\varphi _{m+1}\leq C^{s_{1}}\sigma ^{2s_{2}}\varphi _{0}
\end{equation*}%
where 
\begin{equation*}
s_{1}=\tsum_{j=0}^{m}\sigma ^{-j}\text{ and }s_{2}=\tsum_{j=0}^{m}j\sigma
^{-j}.
\end{equation*}%
In the limit as $m\rightarrow \infty $ we obtain an estimate for $\varphi
_{0}$ in terms of $\left\Vert \overset{\symbol{126}}{u}\right\Vert
_{2,\infty },\left\Vert \overset{\symbol{126}}{u_{x}}\right\Vert _{2,2}$ and 
$\kappa $. The final result is obtained by using further iteration arguments
based on a variant of the fundamental equality (7) to estimate\textit{\ }$%
\left\Vert \overset{\symbol{126}}{u}\right\Vert _{2,\infty }$ and $%
\left\Vert \overset{\symbol{126}}{u_{x}}\right\Vert _{2,2}$\textit{.\ }

In the absence of information about the boundary behavior of a weak solution
of (1) it is possible to use iterative arguments similar to those used to
prove the maximum principle to bound $u$ in a subcylinder in terms of the
norm $\left\Vert u\right\Vert _{2,2}$ over the full cylinder. Let $%
(x^{\prime },t^{\prime })$ be a fixed point in the cylinder $Q.$ Let $R(\rho
)$ denote the open cube in $\mathbf{R}^{n}$ with edge length $\rho $
centered at $x^{\prime }$ and define%
\begin{equation*}
Q(\rho )=R(\rho )\times (t^{\prime }-\rho ^{2},t^{\prime }).
\end{equation*}%
The symbol $\left\Vert \cdot \right\Vert _{p,q,\rho }$ will be used to
denote the $L^{p,q}$ norm of a function over the cylinder $Q(\rho )$\textbf{.%
}

\QTP{Body Math}
\textit{LOCAL\ BOUNDEDNESS}: \textit{Let }$u$\textit{\ be a weak solution of
(1) in }$Q$\textit{\ and suppose that }$Q(3\rho )\subset Q.$\textit{\ Then
almost everywhere in }$Q(\rho )$\textit{\ we have}%
\begin{equation*}
\left\vert u(x,t)\right\vert \leq C(\rho ^{-(n+2)/2}\left\Vert u\right\Vert
_{2,2,3\rho }+\rho ^{\theta }k)
\end{equation*}%
\textit{where }$C$\textit{\ is a constant depending only on }$\rho $\textit{%
\ and the structure of (1) and }%
\begin{equation*}
k=\left\Vert f\right\Vert +\left\Vert g\right\Vert +\left\Vert h\right\Vert .
\end{equation*}%
\textit{In particular, weak solutions of (1) are locally essentially bounded}%
.

\QTP{Body Math}
The Harnack inequality for non-negative harmonic functions gives a bound for
the maximum of a harmonic function over an interior subset of the domain of
definition in terms of the minimum over the same subset. For parabolic
equations there is a similar result except that now the subsets must be
separated by a non-empty time interval as explained in Moser's paper [10].
Moser deals with the linear equation%
\begin{equation*}
u_{t}-\left\{ A_{ij}(x,t)u_{x_{i}}\right\} _{x_{j}}=0
\end{equation*}%
while in [5] the Harnack inequality is proved following the outline of
Moser's proof but applied to the nonlinear equation (1).

\QTP{Body Math}
Let $(x^{\prime },t^{\prime })$ be a fixed point in the basic set $Q$. Let $%
Q(\rho )$ be as defined above and define%
\begin{equation*}
Q^{\ast }(\rho )=R(\rho )\times (t^{\prime }-8\rho ^{2},t^{\prime }-7\rho
^{2}),
\end{equation*}%
that is, $Q^{\ast }(\rho )$ translated downward a distance $7\rho ^{2}$.

\QTP{Body Math}
\textit{HARNACK\ INEQUALITY:\ Let }$u$\textit{\ be a non-negative weak
solution of (1) in }$Q$\textit{. Suppose that }$Q(3\rho )\subset Q$\textit{.
Then}%
\begin{equation*}
\max_{Q^{\ast }(\rho )}u\leq C\min_{Q(\rho )}(u+\rho ^{\theta }k)
\end{equation*}%
\textit{where }$C$\textit{\ is a constant depending only on }$\rho $\textit{%
\ and the structure of (1) and }$k=\left\Vert f\right\Vert +\left\Vert
g\right\Vert +\left\Vert h\right\Vert .$

\QTP{Body Math}
Here $\ \min $ and $\max $ stand for the essential minimum and essential
maximum, both of which are finite thanks to the local boundedness theorem.\
Moreover,C is independent of $\left\Vert f\right\Vert ,\left\Vert
g\right\Vert $ and $\left\Vert h\right\Vert $.

\QTP{Body Math}
An important consequence of the local boundedness theorem and the Harnack
inequality is the uniform H\"{o}lder continuity of solutions of equation
(1). We write $X=(x,t),Y=(y,s),$ etc. to denote points in space-time and
introduce a pseudo-distance according to the definition%
\begin{equation*}
\left\vert X\right\vert ^{2}=\left\{ 
\begin{array}{c}
\max (x_{i}^{2},-t/4)\text{ for }t\leq 0 \\ 
\infty \text{ for }t>0%
\end{array}%
.\right.
\end{equation*}%
Thus the set $\left\vert Y-X\right\vert <\rho $ for fixed $X$ is the cylinder%
\begin{equation*}
\left\vert x_{i}-y_{i}\right\vert <\rho ,t-4\rho ^{2}<s\leq t.
\end{equation*}

\QTP{Body Math}
\textit{H\"{O}LDER\ CONTINUITY}: \textit{Suppose that }$u$\textit{\ is a
weak solution of (1) in }$Q$\textit{. Then }$u$\textit{\ is (essentially) H%
\"{o}lder continuous in }$Q$\textit{. Moreover, if }$\left\vert u\right\vert
\leq L$\textit{\ and }$X^{\prime },Y^{\prime }$\textit{\ are points of }$Q$%
\textit{\ with }$s\leq t$\textit{\ then}%
\begin{equation*}
\left\vert u(Y)-u(X)\right\vert \leq H(L+k)\left( \frac{\left\vert
Y-X\right\vert }{R}\right) ^{\alpha }
\end{equation*}%
\textit{where }$H$\textit{\ and }$\alpha $\textit{\ are positive constants
depending only on the structure of (1), }$k=\left\Vert f\right\Vert
+\left\Vert g\right\Vert +\left\Vert h\right\Vert $\textit{\ and }$R$\textit{%
\ is the pseudo-distance from }$X$\textit{\ to the boundary of }$Q$\textit{\
or }$R=1$\textit{\ if this is smaller.}

\QTP{Body Math}
The proof of the H\"{o}lder continuity of $u$ proceeds by an iteration
argument based on the Harnack inequality applied to the oscillation of $u.$%
Since every weak solution of (1) is locally continuous it is meaningful to
consider the value of a solution at a point. In particular, we can derive a
pointwise Harnack inequality which makes explicit the dependence of the
Harnack constant on the domain.

\QTP{Body Math}
The results so far concern solutions defined in the space-time cylinder $Q$.
They are also valid for solutions defined in the strip $S=\mathbf{R}%
^{n}\times (0,T).$There is a drawback here, however, namely the requirement
that the norms $\left\Vert b\right\Vert $ through $\left\Vert h\right\Vert $
be finite over the entire strip. Thus the results as stated do not apply to
a solution in $S$ when, for example, any of the quantities $b$ through $h$
is constant. To remedy this defect we introduce the set $\mathcal{S}$ of all
cylinders of the form $Q(\sigma )$, with $\sigma =\min (1,\sqrt{T}),$
contained in $S$. A function $w=w(x,t)$ defined on $S$ is said to belong to
the class $L^{p,q}(\mathcal{S})$ if%
\begin{equation*}
\sup \left\Vert w\right\Vert _{p,q}<\infty
\end{equation*}%
where the norms are taken over cylinders in the family $\mathcal{S}$. All of
the preceding results continue to hold for solutions of (1) in $S$, where
the coefficients in the structural inequalities are contained in the
respective classes $L^{p,q}(\mathcal{S})$ rather than $L^{p,q}(Q)$.

\QTP{Body Math}
\textit{POINTWISE HARNACK INEQUALITY}: \textit{Let }$u$\textit{\ be a
non-negative weak solution of (1) in }$S,$\textit{\ where the coefficients }$%
b$\textit{\ through }$h$\textit{\ in the structure inequalities are
contained in the appropriate classes }$L^{p,q}(\mathcal{S)}$\textit{. Then
for all points }$(x,t)$\textit{\ and }$(y,s)$\textit{\ in }$S$\textit{\ with 
}$0<s<t<T$\textit{\ we have}%
\begin{equation*}
u(y,s)+k\leq \{u(x,t)+k\}\exp C\left( \frac{\left\vert x-y\right\vert ^{2}}{%
t-s}+\frac{t}{s}\right) .
\end{equation*}%
\textit{Here }$k=\sup \left\Vert f\right\Vert +\sup \left\Vert g\right\Vert
+\sup \left\Vert h\right\Vert $\textit{\ and }$C$\textit{\ depends only on
the structure of (1) and on }$T$\textit{.}

\QTP{Body Math}
The proof of this result is based on the Harnack inequality and uses the
fact that we can talk about the pointwise values of $u$. Note that the
pointwise Harnack inequality holds for the linear equation (5) even when the
coefficients are constant.

\QTP{Body Math}
The following result which gives a lower bound on the growth of a solution
of (1) in $S$ as $t\searrow 0$ is a generalization of a result due to Nash
[11] for the linear equation\ $u_{t}-\left\{ A_{ij}(x,t)u_{x_{i}}\right\}
_{x_{j}}=0.$\ It is a direct consequence of the pointwise Harnack inequality
and plays an essential role in the Gaussian estimates discussed below.

\textit{LIMIT BEHAVIOR}:\textit{\ Let }$u$\textit{\ be a non-negative weak
solution of (1) in the strip }$S$\textit{, where the coefficients in (2) are
contained in the appropriate classes }$L^{p,q}(\mathcal{S})$\textit{.
Suppose that for some }$\alpha >0$\textit{\ we have}%
\begin{equation*}
\mathcal{M=}\inf_{0<t<T}\tint_{\left\vert x\right\vert ^{2}<\alpha
t}u(x,t)dx>0.
\end{equation*}%
\textit{Then there exist positive constants }$C_{1}$\textit{\ and }$C_{2}$%
\textit{\ such that}%
\begin{equation*}
u(x,t)+k\geq C_{1}t^{-n/2}\exp \left( -C_{2}\left\vert x\right\vert
^{2}/t\right) 
\end{equation*}%
\textit{in }$S$\textit{. Here }$C_{1}$\textit{\ depends only on }$\alpha ,%
\mathcal{M},n,T$\textit{\ and the structure of (1), while }$C_{2}$\textit{\
depends only on }$T$\textit{\ and the structure of (1).}

\QTP{Body Math}
In [1], [2] and [3] the results of [5] are applied to the study of
non-negative solutions of the linear equation (5). Many of the results in
[1], [2] and [3] concern the properties of the weak fundamental solution $%
\Gamma (x,t;\xi ,\tau )$ of the homogeneous equation%
\begin{equation}
u_{t}=\left\{ A_{ij}(x,t)u_{x_{i}}+A_{j}(x,t)u\right\}
_{x_{j}}+B_{j}(x,t)u_{x_{j}}+C(x,t)u.  \tag{8}
\end{equation}%
It is shown\ in [2] that the weak fundamental solution has all the essential
properties of the classical fundamental solution and indeed that they
coincide whenever the classical solution exists. Among the key results
proved in [1] and [2] are the following bounds for $\Gamma $.

\QTP{Body Math}
\textit{GAUSSIAN BOUNDS}:\textit{\ Let }$\Gamma (x,t;\xi ,\tau )$\textit{\
be the weak fundamental solution of a uniformly parabolic equation (8) whose
coefficients are measurable and contained in the appropriate classes }$%
L^{p,q}(\mathcal{S})$. \textit{Then there exist positive constants }$\alpha
_{1},\alpha _{2}$\textit{\ and }$C$\textit{\ depending only on }$T$\textit{\
and the bounds for the coefficients\ such that}%
\begin{equation*}
\mathcal{C}^{-1}g_{1}(x-\xi ;t-\tau )\leq \Gamma (x,t;\xi ,\tau )\leq 
\mathcal{C}g_{2}(x-\xi ;t-\tau )
\end{equation*}%
\textit{for all }$(x,t,\xi ,\tau )\in S\times S$\textit{\ with }$t>\tau $%
\textit{, where }$g_{i}(x,t)$\textit{\ is the fundamental solution of the
heat conduction equation }$\alpha _{i}\Delta u=u_{t}$\textit{\ for }$i=1,2$%
\textit{.}

\QTP{Body Math}
It is noteworthy that these bounds do not require any smoothness assumptions
on the coefficients of (8), unlike similar bounds which had previously
appeared in the literature. The proof of the upper bound given in [1] and
[2]\ is independent of the results in [5]. It is based on the Kolmogorov
identity%
\begin{equation*}
\Gamma (x,t;\xi ,\tau )=\tint_{\mathbf{R}^{n}}\Gamma (x,t;\varsigma ,\eta
)\Gamma (\zeta ,\eta ;\xi ,\tau )d\zeta
\end{equation*}%
and a technical estimate\ for the growth at the center of a ball of a
solution which is initially supported in the exterior of that ball. The the
proof of the lower bound is based on the Limit Behavior Theorem. For the
homogenous equation (8), $k=0$ and the crux of the proof is the estimation
of $\mathcal{M}$, which is carried out using the Pointwise Harnack
Inequality.

\QTP{Body Math}
Fabes and Stroock [7] made a deep study of Nash's methods and showed for the
equation $u_{t}-\left\{ A_{ij}(x,t)u_{x_{i}}\right\} _{x_{j}}=0$ that the
Gaussian bounds could be derived directly using an extension of Nash's ideas
without first proving continuity and a Harnack principle. Indeed they showed
that H\"{o}lder continuity and a Harnack inequality could be derived
directly from the Gaussian bounds. Thus these three results-H\"{o}lder
continuity, Harnack inequality and Gaussian estimates-are connected and are
in some sense equivalent. Norris and Stroock [12] extended the results of
[7] to more general equations.

\QTP{Body Math}
Before continuing it is worth noting that in case of the equation $%
u_{t}-\left\{ A_{ij}(x,t)u_{x_{i}}\right\} _{x_{j}}=0$, if the coefficients $%
A_{ij}$ are independent of $t$ and if $n\geq 3$ then%
\begin{equation*}
\dint\limits_{0}^{\infty }\Gamma (x,t;\xi ,0)dt=G(x.\xi ),
\end{equation*}%
where $G(x,\xi )$ is the fundamental solution of the elliptic equation%
\begin{equation*}
\left\{ A_{ij}(x)u_{x_{i}}\right\} _{x_{j}}=0
\end{equation*}%
As is remarked in [1], in this case the constants in the Gaussian estimate
can be chosen independent of $T$ so we can integrate these estimates over $%
\mathbf{R}^{+}$ to obtain%
\begin{equation*}
K^{-1}\left\vert x-\xi \right\vert ^{2-n}\leq G(x,\xi )\leq K\left\vert
x-\xi \right\vert ^{2-n}.
\end{equation*}%
This result was previously derived directly from potential theoretic
considerations by Littman, Stampacchia and Weinberger [8] and by H. Royden
[13].

\QTP{Body Math}
A generalization of the Widder Representation Theorem for the heat
conduction equation [15] is proved in [3]. This result gives a complete
characterization of a non-negatives solution of (8).

\QTP{Body Math}
\textit{WIDDER REPRESENTATION THEOREM}: \textit{Suppose that equation (8) is
uniformly parabolic with measurable coefficients in the appropriate classes }%
$L^{p,q}(S).$\textit{If }$u$\textit{\ is a non-negative weak solution of (8)
in }$S$\textit{\ then there exist a unique non-negative Borel measure }$\rho 
$\textit{\ on }$R^{n}$\textit{\ such that}%
\begin{equation}
u(x,t)=\tint_{\mathbf{R}^{n}}\Gamma (x,t;\xi ,0)\rho (d\xi ),  \tag{9}
\end{equation}%
\textit{where }$\Gamma $\textit{\ is the weak fundamental solution of (8)
and }$\rho $\textit{\ satisfies}%
\begin{equation}
\tint_{\mathbf{R}^{n}}e^{-\sigma \left\vert x\right\vert ^{2}}\rho
(dx)<\infty  \tag{10}
\end{equation}%
\textit{for some }$\sigma >0.$\textit{Moreover, the measure }$\rho $\textit{%
\ is the initial trace of }$u$\textit{, that is,}%
\begin{equation*}
\lim_{t\searrow 0}\tint_{\mathbf{R}^{n}}u(x,t)\psi (x)dx=\tint_{\mathbf{R}%
^{n}}\psi (x)\rho (dx)
\end{equation*}%
\textit{for all }$\psi \in C(R^{n})$\textit{\ such that }$\left\vert \psi
(x)\right\vert \leq Ke^{-\delta \left\vert x\right\vert ^{2}}$\textit{\ for
some constant }$K$\textit{\ and }$\delta >\sigma .$\textit{Conversely, if }$%
u $\textit{\ is given by (9) with a non-negative Borel measure }$\rho $%
\textit{\ satisfying (10) then }$u$\textit{\ is a non-negative weak solution
of (8) with initial trace }$\rho .$

\QTP{Body Math}
Note that if $\rho $ has a density $\mu $ then%
\begin{equation*}
u(x,t)=\tint_{\mathbf{R}^{n}}\Gamma (x,t;\xi ,0)\mu (\xi )d\xi .
\end{equation*}

\QTP{Body Math}
The existence and uniqueness of the representing measure $\rho $ is proved
in [2] and its characterization (10) is a consequence of the Gaussian lower
bound. The proof that a function given by (9) with a measure satisfying (10)
is a non-negative weak solution depends on an approximation argument and the
Gaussian upper bound.

\QTP{Body Math}
Serrin's seminal work on regularity of solutions of quasilinear elliptic
equations [14] led directly to the extensions to quasilinear parabolic
equations in [4] and [5]. These papers, in turn, gave rise to the detailed
study of non-negative solution of linear parabolic equations described in
[1], [2] and [3].

\QTP{Body Math}
$\bigskip $

\QTP{Body Math}
REFERENCES

\QTP{Body Math}
1. \ D. G. Aronson, \textit{Bounds for the fundamental solution of a
parabolic equation}. Bull. A. M. S., \textbf{73}(1967), 890-896.

\QTP{Body Math}
2. \ D. G. Aronson, \textit{Non-negative solutions to linear parabolic
equations. }Ann. Scuola Norm. Sup. Pisa, \textbf{22}(1968), 607-694.

\QTP{Body Math}
3. \ D. G. Aronson, \textit{Non-negative solutions to linear parabolic
equations: An addendum. }Ann. Scuola Norm. Sup. Pisa, \textbf{25}(1970),
1-10.

\QTP{Body Math}
4. \ D. G. Aronson and James Serrin, \textit{A maximum principle for
non-linear parabolic equations}. \textit{\ }Ann. Scuola Norm. Sup. Pisa, 
\textbf{21}(1967), 291-305-

\QTP{Body Math}
5. \ D. G. Aronson and James Serrin,\textit{\ Local behavior of solutions of
quasilinear parabolic equations.} Arch. Rat. Mech. Anal., \textbf{25}(1967),
81-122.

\QTP{Body Math}
6. \ E. De Giorgi. \textit{Sulla differeziabilit\`{a} e l'analiticit\`{a}
delle estremali degli integrali multipli regulari}. Mem. Acad. Sci. Torino
Cl. Sci. Fis. Mat. Nat., \textbf{3}(1957), 25-43.

\QTP{Body Math}
7. E. B.\ Fabes and D. W. Stroock.\textit{\ A new proof of Moser's parabolic
Harnack inequality using the old ideas of Nash.} Arch. Rat. Mech. Anal., 
\textbf{96}(1986), 327-338.

\QTP{Body Math}
8. \ W. Littman, G. Stampacchia and H. F. Weinberger, \textit{Regular points
for elliptic equations with discontinuous coefficients}. Ann. Scuola Norm.
Sup. Pisa, \textbf{17}(1963), 43-77.

\QTP{Body Math}
9. \ J\"{u}rgen Moser, \textit{On Harnack's theorem for elliptic
differential equations. }Comm. Pure\ App. Math., \textbf{14}(1961),
577-591.\ 

\QTP{Body Math}
10. J\"{u}rgen Moser, \textit{A Harnack's inequality for parabolic
differential equations. }Comm. Pure\ App. Math., \textbf{17}(1964), 101-134

\QTP{Body Math}
11. J. Nash, \textit{Continuity of solutions of parabolic and elliptic
equations.} Amer. J. Math., \textbf{80}(1958), 931-954.

\QTP{Body Math}
12. James Norris and Daniel W. Stroock, \textit{Estimates on the fundamental
solution to heat flow with uniformly elliptic coefficients,}\ Proc. London
Math. Soc., \textbf{62}(1991), 373-402.

\QTP{Body Math}
13. H. L. Royden, \textit{The growth of a fundamental solution of an
elliptic divergence structure equation.}\ Studies in Mathematical Analysis
and Related Topics, Stanford U. Press 1962, \ \ 

\QTP{Body Math}
\ \ 333-340.

\QTP{Body Math}
14. James Serrin,\textit{\ Local behavior of solutions of quasilinear
equations.}\ Acta Math., \textbf{111}(1964), 247-302.\ 

\QTP{Body Math}
15. D. V. Widder,\textit{\ Positive temperature on the infinite rod}. Trans.
Amer. Math. Soc., \textbf{55}(1944),85-95.\ \ \ \ \ \ \ \ \ \ \ \ \ \ \ \ \
\ \ \ \ \ \ \ \ \ \ \ \ \ \ \ \ \ \ \ \ \ \ \ \ \ 

\end{document}